\documentclass[12pt,letter]{amsart}

\usepackage{amsmath}
\usepackage{amscd}
\usepackage{amssymb}
\usepackage[all]{xy}

\newtheorem{theorem}{Theorem}[section]
\newtheorem{theorem/definition}{Theorem/Definition}[section]
\newtheorem{Theorem}{Theorem}
\newtheorem{proposition}{Proposition}[section]

\theoremstyle{remark}

\theoremstyle{definition}

\begin{document}
    \title{A simple proof of Witten Conjecture through localization}
        \author{Yon-Seo Kim}
		\address{Department of Mathematics,University of California at Los Angeles,
				Los Angeles, CA 90095-1555, USA}
        \email{yskim@math.ucla.edu}
		\author{Kefeng Liu}
        \dedicatory{Department of Mathematics, UCLA;
			Center of Math Sciences, Zhejiang University}
		\address{Center of Math Sciences, Zhejiang University, Hangzhou, Zhejiang 310027, China;
				Department of Mathematics,University of California at Los Angeles,
				Los Angeles, CA 90095-1555, USA}
		\email{liu@math.ucla.edu, liu@cms.zju.edu.cn}
        \begin{abstract}
			We obtain a system of relations between Hodge integrals with one $\lambda$-class.
			As an application, we show that its first non-trivial relation implies the Witten's
			Conjecture/Kontsevich Theorem \cite{Wit,Kon}.
        \end{abstract}
    \maketitle

	\section{Introduction}
	In this paper, we obtain an alternate proof of the Witten's Conjecture \cite{Wit}
	which claims that the tautological intersections on the moduli space of stable curves
	$\overline{\mathcal{M}}_{g,n}$ is governed by KdV hierarchy.
	It is first proved by M.Kontsevich \cite{Kon} by constructing combinatorial model for
	the intersection theory of $\overline{\mathcal{M}}_{g,n}$ and interpreting the trivalent
	graph summation by a Feynman diagram expansion for a new matrix integral.
	A.Okounkov-R.Pandharipande \cite{Oko-Pan} and M.Mirzakhani \cite{Mirz} gave different approaches
	through the enumeration	of branched coverings of $\mathbb{P}^{1}$ and the Weil-Petersen volume, respectively.
	Recently, M.Kazarian-S.Lando \cite{Kaz-Lando} obtained an algebro-geometric proof by using the ELSV-formula
	to relate the intersection indices of $\psi$-classes to Hurwitz numbers.\\

	Here we take an approach using virtual functorial localization on the moduli space of relative
	stable morphisms $\overline{\mathcal{M}}_{g}(\mathbb{P}^{1},\mu)$ \cite{Li1}.
	$\overline{\mathcal{M}}_{g}(\mathbb{P}^{1},\mu)$ consists of maps from Riemann surfaces of genus
	$g$ and $n=l(\mu)$ marked points to $\mathbb{P}^{1}$ which has prescribed ramification type
	$\mu$ at $\infty\in\mathbb{P}^{1}$. As the result, we obtain a system of relations between linear
	Hodge integrals. It recursively expresses each linear Hodge integral by lower-dimensional ones.
	The first non-trivial relation of this system is 'cut-and-join relation', and is
	of same recursion type as that of single Hurwitz numbers \cite{Li-Zha-Zhe}.
	Moreover, as we increase the ramification degree, we can extract a relation between absolute
	Gromov-Witten invariants from this relation. And we show this relation implies the following
	recursion relation for the correlation functions of topological gravity \cite{Dijkgraaf}:
	\begin{align*}
		\langle\tilde{\sigma}_{n}\prod_{k\in S}\tilde{\sigma}_{k}\rangle_{g}=&
		\sum_{k\in S}(2k+1)\langle\tilde{\sigma}_{n+k-1}\prod_{l\neq k}\tilde{\sigma}_{l}\rangle_{g}
		+\frac{1}{2}\sum_{a+b=n-2}\langle\tilde{\sigma}_{a}\tilde{\sigma}_{b}\prod_{l\in S}\tilde{\sigma}_{l}\rangle_{g-1}\\
		&+\frac{1}{2}\sum_{S=X\cup Y,a+b=n-2,g_{1}+g_{2}=g}
		\langle\tilde{\sigma}_{a}\prod_{k\in X}\tilde{\sigma}_{k}\rangle_{g_{1}}
		\langle\tilde{\sigma}_{b}\prod_{l\in Y}\tilde{\sigma}_{l}\rangle_{g_{2}}&&\cdots(*)
	\end{align*}
	which is equivalent to the Witten's Conjecture/Kontsevich Theorem.
	This recursion relation (*) is also equivalent to the Virasoro constraints; i.e. (*) can be expressed
	as linear, homogeneous differential equations for the $\tau$-function \cite{Dijkgraaf}
	$$\tau(\tilde{t})=\text{exp}\sum_{g=0}^{\infty}\langle\text{exp}\sum_{n}\tilde{t}_{n}\tilde{\sigma}_{n}\rangle_{g}$$
	$$L_{n}\cdot\tau=0,\qquad (n\geq -1)$$
	where $L_{n}$ denote the differential operators
	\begin{align*}
		L_{-1} &= -\frac{1}{2}\frac{\partial}{\partial\tilde{t}_{0}}+
			\sum_{k=1}^{\infty}(k+\frac{1}{2})\tilde{t}_{k}\frac{\partial}{\partial\tilde{t}_{k-1}}
			+\frac{1}{4}\tilde{t}^{2}_{0}\\
		L_{0} &= -\frac{1}{2}\frac{\partial}{\partial\tilde{t}_{1}}+
			\sum_{k=0}^{\infty}(k+\frac{1}{2})\tilde{t}_{k}\frac{\partial}{\partial\tilde{t}_{k}}+\frac{1}{16}\\
		L_{n} &= -\frac{1}{2}\frac{\partial}{\partial\tilde{t}_{n-1}}+
			\sum_{k=0}^{\infty}(k+\frac{1}{2})\tilde{t}_{k}\frac{\partial}{\partial\tilde{t}_{k+n}}
			+\frac{1}{4}\sum_{i=1}^{n}\frac{\partial^{2}}{\partial\tilde{t}_{i-1}\partial\tilde{t}_{n-i}}
	\end{align*}
	As a remark, it is possible that the general recursion relation obtained from our approach
	implies the Virasoro conjecture for a general non-singular projective variety.

	The rest of this paper is organized as follows:
	In section 2, we recall the recursion formula obtained in \cite{OneLambda} and derive cut-and-join
	relation as its special case.
	In section 3, we prove asymptotic formulas for the coefficients in the cut-and-join relation.
	Then we derive first two relations of the system of relations between linear Hodge integrals,
	and show that the cut-and-join relation implies (*).\\

	\it{* Please refer to \cite{OneLambda} for miscellaneous notations.}\rm

	\section{Recursion Formula}
		The following recursion formula was derived in \cite{OneLambda}.
		\begin{theorem}
			For any partition $\mu$ and $e$ with $\vert e\vert<\vert\mu\vert+l(\mu)-\chi$, we have
			\begin{equation}\label{recursion}
				\big[\lambda^{l(\mu)-\chi}\big]\sum_{\mid\nu\mid=\mid\mu\mid}\Phi^{\bullet}_{\mu,\nu}(-\lambda)
				z_{\nu}\mathcal{D}_{\nu,e}^{\bullet}(\lambda)=0
			\end{equation}
		where the sum is taken over all partitions $\nu$ of the same size as $\mu$.
		\end{theorem}
		Here $\big[\lambda^{a}\big]$ means taking the coefficient of $\lambda^{a}$, and
		$\mathcal{D}^{\bullet}_{\nu,e}$ consists of linear Hodge integrals as follows;
		$$\mathcal{D}_{g,\nu,e}=\frac{1}{l(e)!\,\mid\text{Aut }\nu\mid}
		\Big[\prod_{i=1}^{l(\nu)}\frac{\nu_{i}^{\nu_{i}}}{\nu_{i}!}\Big]
		\int_{\overline{\mathcal{M}}_{g,l(\nu)+l(e)}}
		\frac{\Lambda_{g}^{\vee}(1)\prod_{j=1}^{l(e)}(1-\psi_{j})^{e_{j}}}
		{\prod_{i=1}^{l(\nu)}\big(1-\nu_{i}\psi_{i}\big)}$$
		where $\Lambda^{\vee}_{g}(t)$ is the dual Hodge bundle;
		$$\Lambda^{\vee}_{g}(t)=t^{g}-\lambda_{1}t^{g-1}+\cdots+(-1)^{g}\lambda_{g}$$
		Introduce formal variable $p_{i}$, $q_{j}$ such that 
		$p_{\nu}=p_{\nu_{1}}\times\cdots\times p_{\nu_{l(\nu)}}$,
		$q_{e}=q_{e_{1}}\times\cdots\times q_{e_{l(e)}}$
		, and form a generating series to define $\mathcal{D}^{\bullet}_{\nu,e}$ as follows:
		\begin{align*}
			\mathcal{D}(\lambda,p,q)=&\quad\sum_{\mid\nu\mid\geq 1}\sum_{g\geq 0}\lambda^{2g-2+l(\nu)}p_{\nu}q_{e}\mathcal{D}_{g,\nu}\\
			\mathcal{D}^{\bullet}(\lambda,p,q)=&\quad\text{exp}\big(\mathcal{D}(\lambda,p,q)\big)
			=:\sum_{\mid\nu\mid\geq 0}\lambda^{-\chi+l(\nu)}p_{\nu}q_{e}\mathcal{D}^{\bullet}_{\chi,\nu,e}
			=\sum_{\mid\nu\mid\geq 0}p_{\nu,e}q_{e}\mathcal{D}^{\bullet}_{\nu,e}(\lambda)
		\end{align*}
		The convoluted term $\Phi^{\bullet}_{\mu,\nu}(-\lambda)$ consists of double Hurwitz numbers as follows:
		\begin{align*}
			\Phi^{\bullet}_{\nu,\mu}(\lambda) =&\sum_{\chi}H^{\bullet}_{\chi}(\nu,\mu)
			\frac{\lambda^{-\chi+l(\nu)+\l(\mu)}}{(-\chi+l(\nu)+l(\mu))!}
			\qquad\Phi^{\bullet}(\lambda;p^{0}
			,p^{\infty})=1+\sum_{\nu,\mu}\Phi^{\bullet}_{\nu,\mu}(\lambda)p^{0}_{\nu}p^{\infty}_{\mu}
		\end{align*}
		Here $H^{\bullet}_{\chi}(\nu,\mu)$ is the double Hurwitz number with ramification type
		$\nu$,$\mu$ with Euler characteristic $\chi$.
		The recursion formula (\ref{recursion}) was derived by integrating point-classes over the
		relative moduli space $\overline{\mathcal{M}}_{g}(\mathbb{P}^{1},\mu)$,
		and the 'cut-and-join relation' is only the first term in this much more general formula.
		This can also be seen as follows:
		Denote by $J_{ij}(\mu)$, $C_{i}(\mu)$ for the cut-and-join partitions of $\mu$ \cite{Zho1}
		and consider the following identity obtained by localization method:
		\begin{align*}
			0=\int_{\overline{\mathcal{M}}_{g}(\mathbb{P}^{1},\mu)}\text{Br}^{*}\prod_{k=0}^{r-2}(H-k)
			&=\text{Contribution from the graph that is mapped to }p_{r}\\
			+&\text{Contribution from the graphs that are mapped to }p_{r-1}
		\end{align*}
		It is straightforward to show that preimages of $p_{r}$ and $p_{r-1}$ under the branching 
		morphism $\text{Br}:\overline{\mathcal{M}}_{g}(\mathbb{P}^{1},\mu)\longrightarrow\mathbb{P}^{r}$
		are the unique graph $\Gamma_{r}$ and the 'cut-and-join graphs' of $\Gamma_{r}$, respectively.
		Hence we recover the 'cut-and-join relation' as the restriction of (\ref{recursion})
		to the first two fixed points $\{p_{r},p_{r-1}\}$;
		\begin{equation}\label{cut-and-join}
			r\Gamma_{r}=\sum_{i=1}^{n}\big[
			\sum_{j\neq i}\frac{\mu_{i}+\mu_{j}}{1+\delta^{\mu_{i}}_{\mu_{j}}}\Gamma_{J}^{ij}+
			\sum_{p=1}^{\mu_{i}-1}\frac{p(\mu_{i}-p)}{1+\delta^{p}_{\mu_{i}-p}}\Big(
			\Gamma^{i,p}_{C1}+\sum_{g_{1}+g_{2}=g,\nu_{1}\cup\nu_{2}=\nu}\Gamma^{i,p}_{C2}\Big)
			\big]
		\end{equation}
		where $\Gamma$'s are the contributions from 'cut-and-join' graphs defined as follows;
		\begin{itemize}
			\item{Original graph that is mapped to the branching point $p_{r}$}
			$$\Gamma_{r}=\frac{1}{\vert\text{Aut }\mu\vert}\prod_{i=1}^{n}\frac{\mu_{i}^{\mu_{i}}}{\mu_{i}!}
				\int_{\overline{\mathcal{M}}_{g,n}}\frac{\Lambda_{g}^{\vee}(1)}{\prod(1-\mu_{i}\psi_{i})}$$
			\item{Join graph that is obtained by joining $i$-th and $j$-th marked points:}
			$$\Gamma_{J}^{ij}=\frac{1}{\vert\text{Aut }\eta\vert}\prod_{k=1}^{n-1}\frac{\eta_{k}^{\eta_{k}}}{\eta_{k}!}
				\int_{\overline{\mathcal{M}}_{g,n-1}}\frac{\Lambda_{g}^{\vee}(1)}{\prod(1-\eta_{k}\psi_{k})}
				,\qquad\eta\in J_{ij}(\mu)$$
			\item{Cut graph that is obtained by pinching around the $i$-th marked point:}
			$$\Gamma_{C1}^{i}=\frac{1}{\vert\text{Aut }\nu\vert}\prod_{k=1}^{n+1}\frac{\nu_{k}^{\nu_{k}}}{\nu_{k}!}
				\int_{\overline{\mathcal{M}}_{g-1,n+1}}\frac{\Lambda_{g-1}^{\vee}(1)}{\prod(1-\nu_{k}\psi_{k})}
				,\qquad\nu\in C_{i}(\mu)$$
			\item{Cut graph that is obtained by splitting around the $i$-th marked point:}
			$$\Gamma_{C2}^{i}=\Big[\prod_{k=1}^{n+1}\frac{\nu_{k}^{\nu_{k}}}{\nu_{k}!}\Big]
				\prod_{s=1,2}\frac{1}{\vert\text{Aut }\nu_{s}\vert}
				\int_{\overline{\mathcal{M}}_{g_{s},n_{s}}}
				\frac{\Lambda_{g_{s}}^{\vee}(1)}{\prod(1-\nu_{s,k}\psi_{k})}
				,\qquad\nu\in C_{i}(\mu)$$
		\end{itemize}
		As was mentioned in \cite{LLZ1}, this 'cut-and-join relation' (\ref{cut-and-join}) recovers the ELSV
		formula \cite{ELSV} since this relation is of the same type as the recursion formula for single Hurwitz numbers
		\cite{Li-Zha-Zhe}, hence giving the identification of the graph contributions with single Hurwitz numbers:
		$$H_{g,\mu}=\frac{r!}{\vert\text{Aut }\mu\vert}\Big[\prod_{i=1}^{l(\mu)}
		\frac{\mu_{i}^{\mu_{i}}}{\mu_{i}!}\Big]\int_{\overline{\mathcal{M}}_{g,l(\mu)}}
		\frac{\Lambda_{g}(1)}{\prod_{i=1}^{l(\mu)}(1-\mu_{i}\psi_{i})}$$
		which is the ELSV formula.
		When there's no confusion, we will denote by $\eta=\eta^{ij}$ for the join-partition and $\nu=\nu^{i,p}$ for the cut-partition
		of splitting $\mu_{i}=p+(\mu_{i}-p)$ for some $1\leq p< \mu_{i}$. Also denote by $\nu_{1}$ and $\nu_{2}$ for the splitting
		of cut-partition $\nu$ such that $\nu_{1}\cup\nu_{2}=\nu$ with $p\in\nu_{1},\mu_{i}-p\in\nu_{2}$.
		Note that in the $\Gamma_{C2}$-type contribution, unstable vertices (i.e. $g=0$ and $n$=1,2) are included.
		We can also use any set $\{p_{k_{0}},\cdots,p_{k_{n}}\}$, $n>0$ of fixed points and obtain
		relations between linear Hodge integrals. And these can be applied to derive deeper relations.

	\section{Degree Analysis}
		In this section, we study asymptotic behaviour of the 'cut-and-join relation'
		and obtain a system of relations between linear Hodge integrals.
		The Hodge integral terms in the graph contributions can be expanded as follows:
		\begin{equation}\label{expansion}
			\int_{\overline{\mathcal{M}}_{g,n}}\frac{\Lambda_{g}^{\vee}(1)}{\prod(1-\mu_{i}\psi_{i})}
			=\sum_{k}\prod\mu_{i}^{k_{i}}\int_{\overline{\mathcal{M}}_{g,n}}\prod\psi_{i}^{k_{i}}
			+\text{lower degree terms}
		\end{equation}
		where $\tilde{k}=(k_{1},\cdots,k_{n})$ are multi-indices running over condition $\sum k_{i}=3g-3+n$.
		Hence the top-degree terms consist of Hodge-integral of $\psi$-classes and lower degree terms
		involve $\lambda$-classes. This will give a system of relations between Hodge integrals
		involving one $\lambda$-class. More precisely, integrals will be determined recursively by
		either lower-dimensional or lower-degree $\lambda$-class integrals.
		The following asymptotic formula is crucial in degree analysis.
		\begin{proposition}\label{asymptotic}
			As $n\longrightarrow\infty$, we have for $k,l\geq 0$
			\begin{align*}
				&e^{-n}\sum_{p+q=n}\frac{p^{p+k+1}q^{q+l+1}}{p!q!} &&\longrightarrow
					\frac{1}{2}\Big[\frac{(2k+1)!!(2l+1)!!}{2^{k+l+2}(k+l+2)!}\Big]n^{k+l+2}+o(n^{k+l+2})\\
				&e^{-n}\sum_{p+q=n}\frac{p^{p+k+1}q^{q-1}}{p!q!} &&	\longrightarrow
					\frac{n^{k+\frac{1}{2}}}{\sqrt{2\pi}}-\Big[\frac{(2k+1)!!}{2^{k+1}k!}\Big]n^{k}+o(n^{k})\\
			\end{align*}
			\begin{proof}
				Let $m$ be an integer such that $1<m<n$ and consider three ranges of $p,q$ as follows:
				\begin{align*}
				&R_{l} = \{\,\,(p,q)\,\,\vert\,\, p>n-m\text{ and }q<m\}\\
				&R_{c} = \{\,\,(p,q)\,\,\vert\,\, m\leq p,q \leq n-m\,\,\}\\
				&R_{r} = \{\,\,(p,q)\,\,\vert\,\, p<m\text{ and }q>n-m\}
				\end{align*}
				Recall the Stirling's formula;
				$$n!=\frac{\sqrt{2\pi}n^{n+1/2}}{e^{n}}\Big(1+\frac{1}{12n}+\cdots\Big)$$
			For the summation over $R_{c}$, let $m=n\epsilon$ and $p=nx$
			for some $\epsilon,x\in\mathbb{R}_{>0}$ so that $m,p\in\mathbb{N}$, then we have
			\begin{align*}
				e^{-n}\sum_{p=m}^{n-m}&\frac{p^{p+k+1}}{p!}\frac{q^{q+l+1}}{q!}
				= \sum_{p=m}^{n-m}\frac{1}{2\pi}p^{k+\frac{1}{2}}q^{l+\frac{1}{2}}\big[1+o(1)\big]\\
				&=\frac{n^{k+l+2}}{2\pi}\sum_{p=m}^{n-m}x^{k+\frac{1}{2}}(1-x)^{l+\frac{1}{2}}
					\frac{1}{n}+o(n^{k+l+2})\\
				&\longrightarrow\frac{n^{k+l+2}}{2\pi}\int_{\epsilon}^{1-\epsilon}x^{k+\frac{1}{2}}
					(1-x)^{l+\frac{1}{2}}dx+o(n^{k+l+2})\qquad\text{as }n\text{ goes to }\infty\\
				&=\frac{n^{k+l+2}}{2\pi}\frac{(2k+1)!!(2l+1)!!}{(2(k+l)+3)!!}\int_{\epsilon}^{1-\epsilon}
					\frac{(1-x)^{k+l+\frac{3}{2}}}{\sqrt{x}}dx+o(n^{k+l+2})+O(\sqrt{\epsilon})\\
				&=\frac{1}{2}\Big[\frac{(2k+1)!!(2l+1)!!}{2^{k+l+2}(k+l+2)!}\Big]n^{k+l+2}
					+o(n^{k+l+2})+O(\sqrt{\epsilon})
			\end{align*}
			As $n\longrightarrow\infty$, we can send $\epsilon\longrightarrow 0$.
			For the summation over $R_{l}$ and $R_{r}$, the top-degree terms belong to $O(n^{k+1/2})$ and $O(n^{l+1/2})$, respectively.
			Since we assume $k,l\geq 0$, both cases belong to $o(n^{k+l+2})$, and this proves the first formula.
			For the second formula, $R_{l}$ has highest order of $n^{k+1/2}$ and one can show that the leading term in the
			asymptotic behaviour is $n^{k+1/2}/\sqrt{2\pi}$. After integration by parts,
			$R_{c}$ gives the second highest term in the asymptotic behaviour
			\begin{align*}
				e^{-n}\sum_{p=m}^{n-1}&\frac{p^{p+k+1}}{p!}\frac{q^{q-1}}{q!}
				= \sum_{p=m}^{n-1}\frac{1}{2\pi}p^{k+\frac{1}{2}}q^{l-\frac{3}{2}}\big[1+o(1)\big]
				=\frac{n^{k}}{2\pi}\sum_{p=m}^{n-1}x^{k+\frac{1}{2}}(1-x)^{-3/2}
					\frac{1}{n}+o(n^{k})\\
				&\longrightarrow\frac{n^{k}}{2\pi}\int_{\epsilon}^{1}x^{k+\frac{1}{2}}
					(1-x)^{-3/2}dx+o(n^{k})\qquad\text{as }n\text{ goes to }\infty\\
				&=\frac{n^{k+1/2}}{\sqrt{2\pi}}-\frac{n^{k}}{2\pi}(2k+1)\int_{\epsilon}^{\delta}
					\frac{x^{k-\frac{1}{2}}}{\sqrt{1-x}}dx+o(n^{k})\\
				&=\frac{n^{k+1/2}}{\sqrt{2\pi}}-\Big[\frac{(2k+1)!!}{2^{k+1}k!}\Big]n^{k}
					+o(n^{k})+O(\sqrt{\epsilon})
			\end{align*}
			This proves the second formula.
			\end{proof}
		\end{proposition}
		Let $\mu_{i}=Nx_{i}$ for some $x_{i}\in\mathbb{R}$ and $N\in\mathbb{N}$.
		By taking general values of $x_{i}$, we can assume, without loss of generality, that
		$\vert\text{Aut }\mu\vert=1$.
		As the ramification degree tends to infinity, i.e. as $N\longrightarrow\infty$, the Hodge
		integral expansion (\ref{expansion}) tends to
		\begin{align*}
			\prod_{i=1}^{n}\frac{\mu_{i}^{\mu_{i}+k_{i}}}{\mu_{i}!}\int_{\overline{\mathcal{M}}_{g,n}}
			\prod\psi_{i}^{k_{i}}+O(e^{N}N^{m-1})
			\longrightarrow e^{\vert\mu\vert}\prod_{i=1}^{n}
			\frac{\mu_{i}^{k_{i}-1/2}}{\sqrt{2\pi}}\int_{\overline{\mathcal{M}}_{g,n}}
			\prod\psi_{i}^{k_{i}}+O(e^{N}N^{m-1})
		\end{align*}
		where $m=3g-3+n-(n/2)$ is the highest degree of $N$ in (\ref{expansion}).
		Same expansion applies to each term in (\ref{cut-and-join}).
		By taking out the common factor $e^{\vert\mu\vert}$ and applying the asymptotic formula
		(\ref{asymptotic}), we find that
		\begin{align*}
			r\Gamma_{r} &=
				N^{m+1}\Big[(x_{1}+\cdots+x_{n})\prod_{i=1}^{n}\frac{x_{i}^{k_{i}-1/2}}{\sqrt{2\pi}}\int_{
				\overline{\mathcal{M}}_{g,n}}\prod_{i=1}^{n}\psi_{i}^{k_{i}}\Big]+O(N^{m})\\
			\Gamma_{C1}^{i} &= \frac{N^{m+1/2}}{2}
				\sum_{k+l=k_{i}-2}\frac{(2k+1)!!(2l+1)!!}{2^{k+l+2}(k+l+2)!}
				x_{i}^{k+l+2}\prod_{j\neq i}\frac{x_{j}^{k_{j}-1/2}}{\sqrt{2\pi}}\Big[
				\int_{\overline{\mathcal{M}}_{g-1,n+1}}\psi_{1}^{k}\psi_{2}^{l}\prod_{j\neq i}
				\psi_{j}^{k_{j}} \\
				&+\sum_{g_{1}+g_{2}=g, \nu_{1}\cup\nu_{2}=\nu}\int_{\overline{\mathcal{M}}_{g_{1},
				n_{1}}}\psi_{1}^{k}\prod\psi_{j}^{k_{j}}\int_{\overline{\mathcal{M}}_{g_{2},n_{2}}}
				\psi_{1}^{l}\prod\psi_{j}^{k_{j}}\Big]+O(N^{m}) \\
			\Gamma_{C2}^{i} &=N^{m+1/2}\prod_{j\neq i}\frac{x_{j}^{k_{j}-1/2}}{\sqrt{2\pi}}
				\Big[\sqrt{N}\frac{x_{i}^{k_{i}+1/2}}{\sqrt{2\pi}}\int_{\overline{\mathcal{M}}_{g,n}}
				\prod_{l=1}^{n}\psi_{l}^{k_{l}}-\frac{(2k_{i}+1)!!}{2^{k_{i}+1}k_{i}!}x_{i}^{k_{i}}
				\int_{\overline{\mathcal{M}}_{g,n}}\prod_{l=1}^{n}\psi_{l}^{k_{l}}\Big]+O(N^{m}) \\
			\Gamma_{J}^{ij} &=
				N^{m+1/2}\frac{(x_{i}+x_{j})^{k_{i}+k_{j}-1/2}}{\sqrt{2\pi}}
				\prod_{l\neq i,j}\frac{x_{l}^{k_{l}-1/2}}{\sqrt{2\pi}}\int_{\overline{\mathcal{M}}_{g,n-1}}
				\psi^{k_{i}+k_{j}-1}\prod_{l\neq i,j}\psi_{l}^{k_{l}}+O(N^{m})
		\end{align*}
		Putting them together in the 'cut-and-join relation' (\ref{cut-and-join}) yields
		a system of relations between Hodge integrals with one $\lambda$-class as follows:
		First, we have a system of relations given by the spectrum of $N$-degree. Secondly,
		each relation given by some fixed $N$-degree stratum can be viewed as a polynomial in $x_{i}$'s;
		$$R_{\tilde{m}}(x_{1},\cdots,x_{n}) = \sum_{(s_{1},\cdots,s_{n})}C(s_{1},\cdots,s_{n})
		x_{1}^{s_{1}}\cdots x_{n}^{s_{n}}$$
		where $\tilde{m}$ is a half integer less than or equal to $m+1$ and the coefficient $C(s_{i})$ of
		the homogeneous polynomial $x_{1}^{s_{1}}\cdots x_{n}^{s_{n}}$ involves linear Hodge integrals.
		Since $x_{i}$'s are independent variables, we obtain vanishing relations for each of $C(s_{i})$'s.
		In particular, the first few vanishing relations are given as follows:
		\begin{itemize}
		\item[$\bullet$] For $N^{m+1}$-stratum, we have a trivial identity:
			\begin{align*}
			(x_{1}+\cdots+x_{n})\prod\frac{x_{i}^{k_{i}-1/2}}{\sqrt{2\pi}}\int_{\overline{\mathcal{M}}
			_{g,n}}\prod\psi_{i}^{k_{i}}
			-(x_{1}+\cdots+x_{n})\prod\frac{x_{i}^{k_{i}-1/2}}{\sqrt{2\pi}}
			\int_{\overline{\mathcal{M}}_{g,n}}\prod\psi_{i}^{k_{i}}=0
			\end{align*}
		\item[$\bullet$] From $N^{m+1/2}$-stratum, we obtain a relation between cut-and-join graphs:
			\begin{align*}
				&\sum_{i=1}^{n}\Big[\frac{(2k_{i}+1)!!}{2^{k_{i}+1}k_{i}!}x_{i}^{k_{i}}
				\prod_{j\neq i}\frac{x_{j}^{k_{j}-1/2}}{\sqrt{2\pi}}\int_{\overline{\mathcal{M}}_{g,n}}
				\prod\psi_{j}^{k_{j}}\\
				&-\sum_{j\neq i}\frac{(x_{i}+x_{j})^{k_{i}+k_{j}-1/2}}{\sqrt{2\pi}}\prod_{l\neq i,j}
				\frac{x_{l}^{k_{l}-1/2}}{\sqrt{2\pi}}\int_{\overline{\mathcal{M}}_{g,n-1}}
				\psi^{k_{i}+k_{j}-1}\prod\psi_{l}^{k_{l}}\\
				&-\frac{1}{2}\sum_{k+l=k_{i}-2}\frac{(2k+1)!!(2l+1)!!}{2^{k+l+2}(k+l+2)!}
				x_{i}^{k_{i}}\prod_{j\neq i}\frac{x_{j}^{k_{j}-1/2}}{\sqrt{2\pi}}
				\Big[\int_{\overline{\mathcal{M}}_{g-1,n+1}}\psi_{1}^{k}\psi_{2}^{l}\prod\psi_{j}^{k_{j}}\\
				&+\sum_{g_{1}+g_{2}=g, \nu_{1}\cup\nu_{2}=\nu}\int_{\overline{\mathcal{M}}_{g_{1},n_{1}}}
				\psi_{1}^{k}\prod\psi_{j}^{k_{j}}\int_{\overline{\mathcal{M}}_{g_{2},n_{2}}}\psi_{1}^{l}
				\prod\psi_{j}^{k_{j}}\Big]\Big] = 0\quad\cdots\text{(**)}\\
			\end{align*}
		\item[$\bullet$] Lower degree strata will give relations for Hodge integrals involving
				non-trivial $\lambda$-class in terms of lower-dimensional ones. For example, the
				relation given by the $N^{m}$-stratum recovers the $\lambda_{1}$-expression.
		\end{itemize}
		And the first non-trivial relation (**) implies the Witten's Conjecture (*):
		\begin{Theorem}
			The relation (**) implies (*).
		\begin{proof}
		Introduce formal variables $s_{i}\in\mathbb{R}_{>0}$ and recall the Laplace Transformation:
		$$\int_{0}^{\infty}\frac{x^{k-1/2}}{\sqrt{2\pi}}e^{-x/2s}dx=(2k-1)!!\,\,s^{k+1/2},\qquad
		\int_{0}^{\infty}x^{k}e^{-x/2s}dx = k!\,\,(2s)^{k+1}$$
		Applying Laplace Transformation to the $N^{m+1/2}$-stratum gives the following relation:
		\begin{align*}
			\sum_{i=1}^{n}\Big[&\quad
			s_{i}^{k_{i}+1}(2k_{i}+1)!!\prod_{j\neq i}s_{j}^{k_{j}+1/2}(2k_{j}-1)!!
			\int_{\overline{\mathcal{M}}_{g,n}}\prod\psi_{l}^{k_{l}} \\
			& -\sum_{a+b=k_{i}-2}s_{i}^{k_{i}+1}(2a+1)!!(2b+1)!!\prod_{j\neq i}s_{j}^{k_{j}+1/2}(2k_{j}-1)!!\\
			&\times\big(\int_{\overline{\mathcal{M}}_{g-1,n+1}}\psi_{1}^{a}\psi_{2}^{b}\prod\psi_{l}^{k_{l}}
			+\sum_{g_{1}+g_{2}=g,\cdots}
			\int_{\overline{\mathcal{M}}_{g_{1},n_{1}}}\psi^{a}\prod\psi_{l}^{k_{l}}
			\int_{\overline{\mathcal{M}}_{g_{2},n_{2}}}\psi^{b}\prod\psi_{l}^{k_{l}}\big)\\
			&-\sum_{j\neq i}\frac{(2w+1)!!}{\sqrt{s_{i}}+\sqrt{s_{j}}}
			\big( s_{i}s_{j}^{w+2}+s_{i}^{3/2}s_{j}^{w+3/2}+\cdots+s_{i}^{w+2}s_{j}\big)\\
			&\times\prod_{l\neq i,j} s_{l}^{k_{l}+1/2}(2k_{l}-1)!!
			\int_{\overline{\mathcal{M}}_{g,n-1}}\psi^{w}\prod\psi_{l}^{k_{l}}\quad\Big]=0\\
		\end{align*}
		where $w=k_{i}+k_{j}-1$. The last term is derived from direct integration;
		\begin{align*}
			&\frac{N^{k+\frac{1}{2}}}{\sqrt{2\pi}}\int_{0}^{\infty}\int_{0}^{\infty}(x_{i}+x_{j})^{k+\frac{1}{2}}e^{-x_{i}y_{i}}e^{-x_{j}y_{j}}dx_{i}dx_{j}
			=\frac{N^{k+\frac{1}{2}}}{2\sqrt{2\pi}}\int_{0}^{\infty}\int_{-r}^{r}r^{k+\frac{1}{2}}e^{-\frac{r+s}{2}y_{i}}e^{-\frac{r-s}{2}y_{j}}dsdr\\
			&=\frac{N^{k+\frac{1}{2}}}{2\sqrt{2\pi}}\int_{0}^{\infty}\Big[\int_{-r}^{r}
				e^{\frac{y_{j}-y_{i}}{2}s}ds\Big]r^{k+\frac{1}{2}}e^{-\frac{y_{i}+y_{j}}{2}r}dr
			=\frac{N^{k+\frac{1}{2}}}{\sqrt{y_{i}}+\sqrt{y_{j}}}\frac{(2k+1)!!}{(2y_{i}y_{j})^{k+\frac{3}{2}}}\Big[y_{i}^{k+1}+y_{i}^{k+\frac{1}{2}}y_{j}^{\frac{1}{2}}+\cdots+y_{j}^{k+1}\Big]
		\end{align*}
		under change of variable $r=x_{i}+x_{j}$ and $s=x_{i}-x_{j}$.
		Considering this as a polynomial in $s_{i}$'s, we can isolate out coefficients to obtain
		\begin{align*}
			&(\#)\cdots(2k_{i}+1)!!\prod_{j\neq i}(2k_{j}-1)!!\int_{\overline{\mathcal{M}}_{g,n}}
			\prod\psi_{l}^{k_{l}} = \sum_{j\neq i}(2w+1)!!\prod_{l\neq i,j}(2k_{l}-1)!!
			\int_{\overline{\mathcal{M}}_{g,n-1}}\psi^{w}\prod_{l\neq i,j}\psi_{l}^{k_{l}}+\\
			&\sum_{a+b=k_{i}-2}(2a+1)!!(2b+1)!!\Big[\int_{\overline{\mathcal{M}}_{g-1,n+1}}
			\psi^{a}\psi^{b}\prod_{l\neq i}\psi_{l}^{k_{l}}+\sum\int_{\overline{\mathcal{M}}_{g_{1},n_{1}}}
			\psi^{a}\prod\psi_{l}^{k_{l}}\int_{\overline{\mathcal{M}}_{g_{2},n_{2}}}\psi^{b}\prod\psi_{l}^{k_{l}}\quad\Big]
		\end{align*}
		The reason for getting $1$ as coefficient in the Join-case is due to the following expansion
		\begin{align*}
			&\frac{1}{\sqrt{s_{i}}+\sqrt{s_{j}}}
			(s_{i}s_{j}^{w+2}+s_{i}^{3/2}s_{j}^{w+3/2}+\cdots+s_{i}^{w+2}s_{j}) \\
			&=\frac{1}{\sqrt{s_{j}}}(1-\sqrt{\frac{s_{i}}{s_{j}}}+\frac{s_{i}}{s_{j}}-(\frac{s_{i}}{s_{j}})^{3/2}+\cdots)
			(s_{i}s_{j}^{w+2}+s_{i}^{3/2}s_{j}^{w+3/2}+\cdots+s_{i}^{w+2}s_{j}) \\
			&=\cdots+1\cdot s_{i}^{k_{i}+1}s_{j}^{k_{j}+1/2}+\cdots
		\end{align*}
		In the notations of (*), we have $\tilde{\sigma}_{n}=(2n+1)!!\sigma_{n}=(2n+1)!!\psi^{n}$ and
		$$\langle\tilde{\sigma}_{k_{1}}\cdots\tilde{\sigma}_{k_{n}}\rangle_{g}
		=\Big[\prod_{i=1}^{n}(2k_{i}+1)!!\Big]
		\int_{\overline{\mathcal{M}}_{g,n}}\psi_{1}^{k_{1}}\cdots\psi_{n}^{k_{n}}$$
		After multiplying a common factor $\prod_{l\neq i}(2k_{l}+1)$ on both sides of (\#), we obtain
		\begin{align*}
			\langle\tilde{\sigma}_{n}\prod_{k\in S}\tilde{\sigma}_{k}\rangle_{g}=&
			\sum_{k\in S}(2k+1)\langle\tilde{\sigma}_{n+k-1}\prod_{l\neq k}\tilde{\sigma}_{l}\rangle_{g}
			+\frac{1}{2}\sum_{a+b=n-2}\langle\tilde{\sigma}_{a}\tilde{\sigma}_{b}\prod_{l\in S}\tilde{\sigma}_{l}\rangle_{g-1}\\
			&+\frac{1}{2}\sum_{S=X\cup Y,a+b=n-2,g_{1}+g_{2}=g}
			\langle\tilde{\sigma}_{a}\prod_{k\in X}\tilde{\sigma}_{k}\rangle_{g_{1}}
			\langle\tilde{\sigma}_{b}\prod_{l\in Y}\tilde{\sigma}_{l}\rangle_{g_{2}}
		\end{align*}
		which is the desired recursion relation (*).
		The factor $2k+1$ comes from missing $j$-th marked point in the Join-graph contribution, and the
		extra $1/2$-factor on Cut-graph contributions is due to graph counting conventions.
		Hence we derived Witten's Conjecture / Kontsevich Theorem through localization on the relative moduli space.
		\end{proof}
		\end{Theorem}

\end{document}